# A FURTHER GENERALIZATION OF THE "GALE-NIKAIDO-KUHN-DEBREU" MARKET EQUILIBRIUM THEOREM

RANJIT VOHRA

ABSTRACT. We extend recent generalizations, by Cornet et al. and by Yannelis, of the classical result of Gale, Nikaido, Kuhn & Debreu (the "GNKD theorem) regarding a rigorous proof of the existence of market equilibrium. We do this by now broadening the applicability of their results, which were limited to economies with commodity space that can be modeled by (real) locally convex Hausdorff spaces, to the wider class of economies with commodity spaces describable by any (real) Hausdorff topological vector space with nontrivial continuous dual.

## 1. INTRODUCTION

As observed in Cornet et al. [7], "One of the fundamental problems in economic theory is the existence of a competitive economic equilibrium." A rigorous proof of this existence, loosely called the GNKD approach (honoring pioneering work by Gale [12], Nikaido [21], Kuhn [18], Debreu [10] and others in the mid 1950's) boils down to proving that the 'excess demand correspondence' of the economy contains 0 at some distinguished price vector $p^*$ (called an *equilibrium price*, or market − clearing price, of the economy). For elaboration of this idea see Ch.18 & 20 in [5]. Modern proofs of GNKD − type theorems, including important contributions by [7], [16] and [28], greatly extended the original Euclidean − space framework of the classical GNKD theorems by proving the equilibrium results for a wider class of (possibly infinite dimensional) commodity spaces, namely the real Hausdorff locally convex topological vector spaces.

The present paper further generalizes previous results by expanding the applicability of the GNKD theorem to any real Hausdorff topological vector space that has a nontrivial continuous dual. This class of real Hausdorff TVSs includes all real Hausdorff locally convex spaces, of course, but also many non − locally convex spaces, such as all sequence spaces $\ell^p$ for $0 < p < 1$ (these have continuous dual isomorphic to $\ell^\infty$); the function space $C(\mathrm{S}^1)$ of real functions on the unit circle, equipped with the topology of convergence in measure, and the Hardy spaces $H^p$ for $0 < p < 1$. We exclude any TVS X with continuous dual $\mathrm{X'} = \{0\}$, e.g. the Lebesgue spaces $L^p([0,1])$ for $0 < p < 1$, because such pathological spaces are inconsistent with the Hahn-Banach extension theorem (used in our Proposition below), and they make the Alaoglu theorem useless for our purposes (see Claim 1 of our Theorem below). The upshot is that although our analysis (like that of [28] & [7]) rests upon the theory of dual systems of the form $(\mathrm{X}, \mathrm{X}^{\#})$, we let (*i*) $\mathrm{X} = (\mathrm{X}, \tau)$ be an arbitrary Hausdorff TVS; (*ii*) the continuous dual of X, namely $\mathrm{X}^{\#} = (\mathrm{X}^{\#}, \sigma(\mathrm{X}^{\#}, \mathrm{X}))$ be the nontrivial continuous dual of X; where $X^{\#}$ is equipped with the weak-* topology $\sigma(X^{\#}, X)$; and (*iii*) the bilinear function of the system given by $(x, \ell) \mapsto \ell(x) \in \mathbb{R}$. See Ch.13, 16 in [1] for details.

## 2. REMINDERS

A **correspondence** $\varphi$ from a set X to a set Y, denoted $\varphi : X \twoheadrightarrow Y$, is a rule that assigns to each $x$ in X a (possibly empty) subset $\varphi(x)$ of Y. Suppose X is a vector space and A is a nonempty subset of X. Then a correspondence $\psi : A \twoheadrightarrow X$ is a **KKM correspondence** on A if for every finite subset $\{x_1, ..., x_n\}$ of A, $\text{co}\{x_1, ..., x_n\} \subseteq \bigcup_{i=1}^{n} \psi(x_i)$. See p.577*ff*, [1] for elaboration. We recall the following:

FAN'S KKM THEOREM (K. Fan, [11]): Let A be a subset of a Hausdorff TVS X, and let $\psi : A \twoheadrightarrow X$ be a KKM correspondence. If $\psi$ is closed-valued and if $\psi(x)$ is compact for some $x \in X$, then $\bigcap_{x \in F} \Gamma(x)$ is nonempty and compact. *Proof*: See p.578, [1].

Suppose X is a Hausdorff TVS. Then the following are equivalent: strict separation in X of disjoint nonempty closed convex sets, one of which has an interior point; and strict separation in X of a closed convex set from a point outside that set.

Given a dual system $\langle X, X^{\#} \rangle$, the **polar** of any nonempty subset A of X, denoted $A^0$, is the subset of the dual space $X^{\#}$ defined by $A^O := \{\ell \in X^{\#} : |\ell.a)| \leq 1, \forall a \in A\}$. Any polar contains 0, and is always convex, balanced and weak-* closed. If $A \subseteq B$ then $A^0 \supseteq B^0$; see p.215*ff*, [1] & Wikipedia, *Polar set*.

ALAOGLU'S THEOREM: Let X be a TVS (not necessarily locally convex or Hausdorff) with continuous dual space $X^{\#}$. Then the polar of any neighborhood U of the origin in X, namely $U^0 := \{\ell \in X^{\#} : |\ell \cdot u)| \leq 1, \forall u \in U\}$, is compact in the weak-* topology. See p.236, [20].

Weak and weak-* topologies are locally convex. The algebraic dual X* of any vector space X **separates points** of X (or is **total**). Any total subspace $X^{\#}$ of X* is closed in the weak -* topology $\sigma(X^*, X)$ on X*. If $X^{\#}$ is total, then the weak topology $\sigma(X, X^{\#})$ on X is Hausdorff.

"PROPER" vs. "STRICT" SEPARATION BY HYPERPLANE. Let A, B be subsets of a vector space X. We say[1] (see p.198, [A+Bo]) that the hyperplane $H = [\ell = \alpha]$

- **separates** A & B if $A \subseteq [\ell \leq \alpha]$ & $B \subseteq [\ell \geq \alpha]$ or vice-versa.
- **strictly separates** A & B if (i) H separates A & B; and (ii) $A \subset [\ell < \alpha]$ & $B \subset [\ell > \alpha]$ or vice-versa.

If the continuous dual X' of a real TVS X is nontrivial, i.e. if $X' \neq \{0\}$, then $|X'|$ is uncountable.

[1] CAUTION: The terminology is not universal; we use definitions & notation in p.198*ff*, [1].

## 3. THE MAIN RESULT

We will need the following preliminary results:

PROPOSITION [STRICT SEPARATION]. Suppose:

- $X = (X, \tau)$ is a real Hausdorff TVS;
- $X^{\#}$, the continuous dual of X, is not trivial, i.e, vector space $X^{\#} \neq \{0\}$, and is equipped with the weak-* topology $\sigma(X^{\#}, X)$;
- the duality function $(X, X^{\#}) \to \mathbb{R}$ is given by $(x, p) \mapsto p \cdot x$;
- $A \subset X$ is closed, convex and has nonempty interior; &
- $x_0 \in X \backslash A$

Then there exists a strict separation of A and $\{x_0\}$.

*Proof.* We may assume, without loss of generality, that $0 \in \text{Int}(A)$, as this can be achieved, if necessary, by a suitable translation. Since A is convex and absorbing, the *Minkowski gauge* $\rho_A$ on X defined by A is the map $\rho_A : X \to [0, \infty)$, $x \mapsto \inf\{\lambda \in \mathbb{R} : x \in \lambda A\}$. Evidently $\rho_A$ is sublinear (Lemma 5.50 (2) in [1]). Furthermore, since $\text{Int}(A) = \{x \in X : \rho_A(x) < 1\}$, the map $\rho_A$ is bounded above on a neighborhood of 0, hence is continuous (Th. 7.7.2(c) in [20)]. Note also that $A = \{x \in X : \rho_A(x) \leq 1\}$, and – since $x_0 \in X \backslash A$ – we have $\rho_A(x_0) > 1$. Consider the 1 – dimensional subspace M of X that is spanned by $x_0$, and define a linear functional $\ell : M \to \mathbb{R}$ by $\ell(tx_0) := t \cdot \rho_A(x_0)$. One checks that the relation $\ell(m) \leq \rho_A(m)$ holds for each $m \in M$, so we know by the *Hahn-Banach Extension Theorem* that there is a linear extension $L : X \to \mathbb{R}$ of $\ell$ such that

(i) $L(x_0) = \ell(x_0) = \rho_A(x_0) > 1$; &

(ii) $L(x) \leq \rho_A(x)$ holds for all $x \in X$.

It follows from (ii) and the continuity of $\rho_A$ that $L$ is bounded above, hence is continuous, i.e. $L \in X^{\#}$. Pick a real number $\gamma \in (1, L(x_0))$, and note that the desired strict separation of A and $\{x_0\}$ then follows by the observation that $L(x_0) > \gamma > L(a)$ holds for all $a \in A$. □

COROLLARY. Let $C \neq X$ be a proper closed convex cone of X, with vertex 0 and with nonempty interior, i.e. $\text{Int}(C) \neq \varnothing$. If $D \subset X$ is nonempty, compact, convex and disjoint from C, then there exists a strict separation of C and D.

*Proof.* Let (C − D) denote the Minkowski difference of C and D. Note that (C − D) is a closed convex subset of X, has nonempty interior, and is disjoint from {0}. Applying the Proposition above to the set (C − D) and 0, we see there is a strict separation of (C − D) and {0}, and this strict separation is equivalent to the desired strict separation of C and D. ☐

The Theorem below, our main result, generalizes Theorem 3.1 in [7], as well as Theorem 3.1 in [28] in that we allow the ambient space to be any "real Hausdorff space X with a non-trivial continuous dual", instead of requiring a (real) Hausdorff locally convex space.

THEOREM. Suppose:

- $\mathrm{X} = (\mathrm{X}, \tau)$ is a real Hausdorff TVS.
- $\mathrm{X}^{\#} \neq \{0\}$ is the continuous dual of X, and $\mathrm{X}^{\#}$ is equipped with the (locally convex) weak-* topology $\sigma(\mathrm{X}^{\#}, \mathrm{X})$, and the duality function $(\mathrm{X}, \mathrm{X}^{\#}) \to \mathbb{R}$ is given by $(x, p) \mapsto p \cdot x$. See p.211 in [1].
- $C \neq \mathrm{X}$ is a proper closed convex cone of X, with vertex 0;
- $v \in \mathrm{Int}(C)$ is an interior point of $C$.
- $C^0 := \{p \in \mathrm{X}^{\#} : p \cdot c \leq 0,\ \forall c \in C\}$ is the **polar cone** of $C$ assumed to be non-degenerate; (in particular $C^0$ is convex, weak-* closed & contains 0);
- $\Delta = \{p \in C^0 : p \cdot v = -1\}$;
- $\xi : \Delta \twoheadrightarrow \mathrm{X}$ is a correspondence (called the "excess demand correspondence") so that
  1. $\xi$ is *upper demicontinuous* with respect to the weak-* topology of $\Delta$ & the topology $\tau$ of X;
  2. $\forall p \in \Delta$, $\xi(p)$ is nonempty, compact and convex;
  3. (Weak Walras' Law) $\forall p \in \Delta$ there exists some $z \in \xi(p)$ such that $p \cdot z \leq 0$.

Then: $\exists p^* \in \Delta$ such that $\xi(p^*) \cap C \neq \varnothing$.

*Proof.* We prove the theorem by establishing three claims.

<u>CLAIM 1</u>. $\Delta$ is nonempty, covex and weak-* compact.

For showing non-emptiness of $\Delta$, we first prove that the polar cone $C^0$, which contains $\Delta$, is nontrivial. Using that $C$ is a nonempty proper subset of X, pick a point $x_0 \in \mathrm{X} \backslash C$. Then the Proposition above, applied to the closed convex set $C$ that has nonempty interior, and to the point $x_0 \in \mathrm{X} \backslash C$, means that there exists a closed hyperplane $[\ell = \alpha]$, where $\ell \in \mathrm{X}^{\#}$ is a nonzero linear functional on X & $\alpha$ is a scalar, such that $\ell \cdot c < \alpha$ holds for all $c \in C$, and $\ell \cdot x_0 > \alpha$. In particular this means $\ell$ is a nonzero element of $C^0$, i.e. $\varnothing \neq C^0 \neq \{0\}$. Using that $C$ is a cone containing the origin 0, the scalar $\alpha$ can be adjusted so that $\ell \cdot c \leq 0$ for all $c \in C$. It will now suffice to prove that $(\ell / - \ell.v) \in \Delta$, where $v \in \mathrm{Int}(C)$ by theorem hypothesis. Since $\ell \neq 0$, we can pick a point $j \in \mathrm{X}$ such that $\ell \cdot j > 0$. Then for small enough $t > 0$, we have $(v + tj) \in C$, so $\ell \cdot (v + tj) \leq 0$ i.e. $\ell \cdot v \leq -t\ell \cdot j < 0$. Hence $(-\ell / \ell \cdot v) \in \Delta$, proving that $\Delta$ is nonempty. Convexity of $\Delta$ is clear from direct calculation.

Next, we prove that $\Delta$ is weak-* compact in $\mathrm{X}^{\#}$. We know $C^0$ is convex, weak-* closed in $\mathrm{X}^{\#}$ & contains 0. Fix a point $v \in \mathrm{Int}(C)$, and pick a balanced $\tau$ - open neighborhood U of the origin 0 small enough that $(v + \mathrm{U}) \subseteq C$. The strategy is to show that $\Delta$ is a closed set contained in a compact subset of $\mathrm{X}^{\#}$. By *Alaoglu's theorem* we know the polar of U, namely the set $\mathrm{U}^0 := \{p \in \mathrm{X}^{\#} : p \cdot u \subseteq [-1,1],\ \forall u \in \mathrm{U}\}$ is weak-* compact in $\mathrm{X}^{\#}$. Since it is clear from its definition that $\Delta$ is weak-* closed in $C^0$ we know (using that $C^0$ is itself weak-* closed in $\mathrm{X}^{\#}$) that $\Delta$ is also weak-* closed in $\mathrm{X}^{\#}$. It remains only to check the inclusion $\Delta \subseteq \mathrm{U}^0$. Pick some $\overline{p} \in \Delta$; this is possible because we have shown above that $\Delta$ is nonempty. We will first show that $\overline{p} \cdot u \subseteq [-1,1]$. From $\overline{p} \in \Delta$ we know $\overline{p} \in C^0 \setminus \{0\}$, such that $\overline{p} \cdot c \leq 0,\ \forall c \in C$. In particular, since $(v + u) \in (v + \mathrm{U}) \subseteq C$ holds for each $u \in \mathrm{U}$, we obtain that $\overline{p} \cdot (v + u) \leq 0$, which says $\overline{p} \cdot u \leq -\overline{p} \cdot v = 1$, i.e. $\overline{p} \cdot u \leq 1$, $\forall u \in \mathrm{U}$. Furthermore, since U is balanced, we may replace $u$ by $-u$ to obtain that $-1 \leq \overline{p} \cdot u$ also holds $\forall u \in \mathrm{U}$. Thus we have the desired inclusion $\overline{p} \cdot u \subseteq [-1,1]$, $\forall u \in \mathrm{U}$. Thus $\Delta \subseteq \mathrm{U}^0$, showing that $\Delta$ is weak-* compact. □

***Remark***. Since $\Delta$ is nonempty, the inclusion $\varnothing \neq \Delta \subseteq \mathrm{U}^0$ proved in Claim 1 shows that the polar $\mathrm{U}^0$, as defined in the proof above, is nontrivial: $\mathrm{U}^0 \neq \{0\}$.

Next, we define a (possibly empty – valued) self – correspondence

$$F : \Delta \twoheadrightarrow \Delta \text{ by}$$

$$p \mapsto \{q \in \Delta : q \cdot x > 0, \forall x \in \xi(p)\}. \qquad (**)$$

A direct calculation shows that $F$ is convex – valued. Furthermore – by the Weak Walras' Law hypotheses – we have that $\forall p \in \Delta$, $p \notin F(p)$. It then follows from Lemma 17.47 in [1] that the *inverse complement correspondence* of $F$, namely

$$\Gamma : \Delta \twoheadrightarrow \Delta$$

$$x \mapsto \Delta \backslash [F^{-1}(x)] = \{p \in \Delta : \exists z \in \xi(p) \text{ with } x \cdot z \leq 0\}$$

is a KKM correspondence.

<u>CLAIM 2</u>. $\exists p^* \in \Delta$ such that $p^* \in \bigcap_{x \in \Delta} \Gamma(x)$.

We will use Fan's KKM theorem (see Reminders, above) to establish that $\bigcap_{x\in\Delta}\Gamma(x)\neq\varnothing$. For this it suffices to check that $\Gamma:\Delta\twoheadrightarrow\Delta$ is weak-* closed – valued, and thus is compact – valued. Pick an arbitrary $q\in\Delta$, and note that $V_q := \{x\in \mathrm{X} : q\cdot x>0\}$ is an open half space in X determined by $q$. Since $\xi:\Delta\twoheadrightarrow \mathrm{X}$ is upper demicontinuous by hypothesis, it follows that the set $\{p\in\Delta:\ \xi(p)\subseteq V_q\}$ is weak-* open in $\Delta$. Evidently $\{p\in\Delta:\ \xi(p)\subseteq V_q\} = (\Delta\backslash\Gamma(q))$, so $\Delta\backslash\Gamma(q)$ is weak -* open in $\Delta$. Therefore $\Gamma(q)$ is weak-* closed in $\Delta$ for each $q$ in $\Delta$, and it then follows from Fan's KKM theorem that $\bigcap_{x\in\Delta}\Gamma(x)\neq\varnothing$. Pick $p^*\in\bigcap_{x\in\Delta}\Gamma(x)$.

<u>CLAIM 3</u>. $\xi(p^*)\cap C\neq\varnothing$.

The claim is proven by contradiction, using the price $p^*$ from claim 2. Suppose, toward a contradiction, that the closed cone $C$ (which is convex & has nonempty interior) and the nonempty convex compact set $\xi(p^*)$ are disjoint: $\xi(p^*)\cap C=\varnothing$. Using the Corollary to the Proposition, above) it follows that $C$ and $\xi(p^*)$ can be strictly separated by a hyperplane $[q^*=\alpha]$ for suitable nonzero $q^*\in \mathrm{X}^{\#}$ and some $\alpha\in\mathbb{R}$. Without loss of generality, we may suppose that $M:=\sup_{c\in C}(q^*\cdot c) < \inf_{z\in\xi(p^*)}(q^*\cdot z) =: N$

Then $M=0$, since $C$ is a cone with vertex 0, so we obtain

$$0=M<N \qquad (***)$$

It follows from the definition of $C^0$ that $q^*\in C^0$, because $\sup_{c\in C}(q^*\cdot c)\leq 0$. Recalling that $q^*\neq 0$ and $q^*\cdot v<0$, put $\lambda:=-(q^*\cdot v)$, so $\lambda>0$. Then, from claim 2, we see that for $(q^*/\lambda)\in\Delta$, $\exists z^*\in\xi(p^*)$ such that $(q^*/\lambda)\cdot z^*\leq 0$, i.e. that $(q^*\cdot z^*)\leq 0$. But then $N=\inf_{z\in\xi(p^*)}(q^*\cdot z)\leq(q^*\cdot z^*)\leq 0$, contradicting (***). We conclude that $\xi(p^*)\cap C\neq\varnothing$ after all. $\square$

## 4. CONCLUDING REMARK

It was observed by Yannelis in [27] that Fan's KKM theorem and the Browder fixed point theorem are essentially equivalent in applications. Owing to this, the proof of Claim 2 in our main result can be reformulated so as to replace Fan's KKM theorem by the following theorem, which is a restatement of Browder's fixed point theorem in [6].

THEOREM (Browder, 1968). Suppose X is a Hausdorff TVS and Y is a nonempty compact convex subset of X. Let $\phi:\mathrm{Y}\twoheadrightarrow\mathrm{Y}$ be a convex – valued correspondence with open lower sections such that for each $y\in\mathrm{Y}$, $y\notin\phi(y)$. Then there is a $y^*\in\mathrm{Y}$ such that $\phi(y^*)=\varnothing$.

Recall the correspondence $F : \Delta \twoheadrightarrow \Delta$, $p \mapsto \{q \in \Delta : q \cdot x > 0, \forall x \in \xi(p)\}$ defined previously (see (**) in preceding section). Observe that Claim 2 in the proof of our Theorem above, namely the assertion "$\exists p^* \in \Delta$ such that $p^* \in \bigcap_{x \in \Delta} \Gamma(x)$" is equivalent to the assertion "$F(p^*) = \varnothing$". We will now give an alternate proof of Claim 2 by proving $F(p^*) = \varnothing$, using Browder's theorem rather than Fan's KKM theorem. Checking the hypotheses of Browder's theorem (as stated above) in our setup, we take $\Delta$ to be the nonempty convex weak-* compact subset of our ambient Hausdorff TVS X, and we note that $(i)$ $F$ is convex − valued, by direct calculation; $(ii)$ $p \notin F(p)$ holds for each $p \in \Delta$, else we contradict the Weak Walras Law of the theorem hypotheses; & $(iii)$ denoting by $V_q$ the open half space $\{z \in \mathrm{X} : q \cdot z > 0\}$, and using that $\xi$ is upper − demicontinuous by hypothesis, the set $F^{-1}(q) = \{p \in \Delta : \xi(p) \subseteq V_q\}$ is open in $\Delta$ for each $q \in \Delta$, meaning that $F$ has open lower sections. Thus all hypotheses of the Browder theorem, as stated above, are satisfied in our setup by the correspondence $F : \Delta \twoheadrightarrow \Delta$, whence we conclude that $\exists p^* \in \Delta$ such that $F(p^*) = \varnothing$.

RANJIT VOHRA
P.O. Box 240, Storrs CT 06268, USA
rvohra@bridgew.edu